\documentclass[reqno]{amsart}
\usepackage{amsmath,amsthm,amssymb,hyperref,graphicx,mathrsfs,mathtools} 

\newtheorem{theorem}{Theorem}
\newtheorem{corollary}{Corollary}
\newtheorem{lemma}{Lemma}
\newtheorem{thmx}{Theorem}
\newtheorem{remark}{Remark}

\allowdisplaybreaks

\begin{document}
\title{Location of the Zeros of Quaternionic Polynomials using Matrix Tools}
\author{N.A.Rather$^1$}
\author{Naseer Ahmad Wani $^2$}
\author{Ishfaq Dar  $^{3*}$}

\address{$^{1,2}$Department of Mathematics, University of Kashmir, Srinagar-190006, India}
\address{$^{3*}$Department of Applied Science, Institute of Technology, University of Kashmir, Srinagar-190006, India}
\email{$^1$dr.narather@gmail.com, $^2$waninaseer570@gmail.com, $^{3*}$ishfaq619@gmail.com,}

\begin{abstract}
Using a variety of matrix techniques, the problem of locating the left  eigenvalues of the  quaternion companion matrices  are investigated in this paper. In a recent paper, Dar et al. \cite{RD}, proved that the zeros of a quaternionic polynomial and the left eigenvalues of corresponding companion matrix are same. In view of this, we use various newly developed matrix techniques to prove various results concerning the location of the zeros of regular polynomials of a quaternionic variable with quaternionic coefficients, which include an extension of the result of A. L. Cauchy as well.
\end{abstract}

\maketitle

\section{\textbf{INTRODUCTION }}
Quaternions are a type of mathematical object that extend the concept of complex numbers into a four-dimensional space. They were first introduced by Irish mathematician William Rowan Hamilton \cite{aa} in the mid-19th century, which are denoted by $\mathbb{H}=\{a=a_{0}+a_{1}i+a_{2}j+a_{3}k:a_{0},a_{1},a_{2},a_{3}\in\mathbb{R}\
 \} $ and $i, j, k$ satisfy $i^2=j^2=k^2=ijk=-1, ~ ij=-ji=k,~ jk=-kj=i,~ ki=-ik=j.$. The set of quaternions is a skew-field and because of non-commutative nature, they differ from complex numbers $\mathbb{C}$ . A   quaternionic number  is denoted by q where $q = {\alpha}+{\beta}i+{\gamma}j+{\delta}k \in\mathbb{H}$, contains one real part ${\alpha}$ and three imaginary parts ${\beta}$, ${\gamma}$ and ${\delta}$. The conjugate of $q$, denoted by $\bar{q}$ is a quaternion $q=\alpha-\beta i-\gamma j-\delta k$ and the norm of $q$ is $|q|=\sqrt{q\bar{q}}=\sqrt{{\alpha}^2+{\beta}^2+{\gamma}^2+{\delta}^2}$. The inverse of each non zero element $q$ of $\mathbb{H}$ is given by $q^{-1}=|q|^{-2}\bar{q}$. We also define the ball $\mathbb{B}(0, r) = \{q \in \mathbb{H} : |q| < r\}$, for $r>0$.
\\Hamilton probably didn't realise, though, that his discovery, quaternions, would eventually be utilised to many fields as in physics,where in  quaternions are correlated to the nature of the universe at the level of quantum mechanics. They lead to elegant expressions of the Lorentz transformations, which form the basis of the modern theory of relativity. In signal processing, Quaternion Fourier Transform (QFT) is a powerful tool. The QFT restores the lost commutative property at the cost of no longer being a division algebra. It can be used, for instance, to embed a watermark in a colour image \cite{dd}. Quaternions are very useful for analyzing situations where rotations in $\mathbb{R}^{3}$ are involved \cite{cc}. Other applications of QFT include face recognition (jointly with Quaternion Wavelet Transform) and voice recognition. Quaternions are very efficient for analyzing situations where rotations in $\mathbb{R}^{3}$ are involved.
In recent years, a lot of research activity has been in progress in the study of the quaternion-related mathematical objects; each year, several research articles are being published in a
wide range of journals and different methodologies are used for various objectives.
The main challenge in the field, currently is the non commutative multiplication of quaternions. However, in
this article, we will try to prove various results concerned with (finite) quaternionic companion matrices,
which in turn yield various results concerning the location of zeros of quaternion polynomials.\\

\maketitle

\section{\textbf{Preliminary}}
\parindent=0mm\vspace{0.in}
\textbf{Quaternion Polynomials:}

A quaternion polynomial of degree n is an expression of the form \\
$f(q)=a_{0}+a_{1}q+...+a_{n}q^n, a_{n}\neq 0$ and $q, a_{i}\in\mathbb{H},~~i=0,...,n.$\\
The addition and the multiplication of two polynomials are defined in the same way as those in the commutative case:\\

 $$(a_{0}+a_{1}q+...+a_{m}q^m)+(b_{0}+b_{1}q+...+b_{n}q^{n})=(a_{0}+b_{0})+(a_{1}+b_{1})q+...+(a_{n}+b_{n})q^{n} $$ \\
 and  \\
 
 $$(a_{0}+a_{1}q+...+a_{m}q^m)(b_{0}+b_{1}q+...+b_{n}q^{n})=c_{0}+c_{1}q+...+c_{m+n}q^{m+n}$$\\
 where $c_{\gamma}=\sum_{\alpha+\beta=\gamma} a_{\alpha}b_{\beta}$.
So, in view of the definition of the multiplication, the indeterminate $q$ is assumed to commute with quaternion coefficients.
 Accordingly it is routine to check that the totality of quaternion polynomials  denoted by $\mathbb{H}$ forms a domain (a ring without zero divisor) \cite{DR}.\\
 
However we cannot define the evaluation of a quaternion polynomial in the same
simple way as that for the commutative case because the direct substitution of $q$ by a
quaternion always yields distinct results for different forms of the polynomial, for instance,
\par Let $f(q)=q^{2}-(i+j)q+k=(q-i)(q-j)$. At $q=i,$
$$ i^{2}-(i+j)i+k=2k\neq 0=(i-i)(i-j)$$
For the sake of integrity, the way of the evaluation must be specified for some certain definite
form of polynomials. There are two main ways of evaluations: the left evaluation and
the right evaluation \cite{bb}.

\parindent=0mm\vspace{0.in}
In this paper, we shall adopt the following definition.
\\

\parindent=0mm\vspace{0.in}
{\bf{Definition:}}
 A quaternionic polynomial of degree n is defined as;\\
 $f(q) = a_nq^{n} +a_{n-1}q^{n-1}...+a_1q + a_0$  or $g(q) = q^{n}a_n +q^{n-1}a_{n-1}...+qa_1 + a_0, $ $a_n \neq 0$ where $ a_i ,q\in \mathbb{H}, i = 1,...,n.$\\
A quaternionic polynomial $p(q)$ is said to be Lacunary  type polynomial if its coefficients skip certain values or are zero at regularly, spaced intervals. Mathematically  given a non-negative integer  $r$ , a polynomial $p(q) = a_nq^{n} +a_{r}q^{r}...+a_1q + a_0,$ or $h(q) = q^{n}a_n +q^{r}a_{r}...+qa_1 + a_0; a_i, q \in \mathbb{H}$, where $a_r \neq 0$,and $1\leq r\leq n-1$ is said to be lacunary type polynomial of degree $n\geq2$ with quaternionic coefficients and $q$ be the quaternionic variable. As usual, if we take the leading coefficient  $a_{n} = 1$, then the  above polynomials are  said to be monic. Without loss of generality, we will take monic polynomials throughout this article.

\parindent=0mm\vspace{0.1in}
{\bf{Quaternionic polynomial matrices:}}As usual, $\mathbb{H}^{r\times s}$ and $\mathbb{H}^{r\times s}[q]$ will, respectively, denote the set of the $r\times s $ matrices with entries in $\mathbb{H}$ and in $\mathbb{H}[q]$. As for polynomials, for simplicity, we may also omit the indeterminate q and write $\textbf{R}$ $\in\mathbb{H}^{r\times r}[q]$.

\parindent=0mm\vspace{0.2in}
{\bf{Quaternion Companion Matrix:}}

The $n \times n$ companion matrix of a monic quaternion polynomial of the form $f(q)=q^n+a_{n-1}q^{n-1}+...+a_1q+a_0$ is given by
$$C_f=\begin{bmatrix}
0& 1& 0& ... & 0& 0\\
0& 0& 1& ... & 0& 0\\
.& .& .& ...& .& .\\
0& 0& 0& ... & 0& 1\\
-a_{0}& -a_{1}& -a_{2}& ...& -a_{n-2}& -a_{n-1}\\
\end{bmatrix}$$

\parindent=0mm\vspace{0.in}Whereas the $n\times n$ companion matrix for a monic quaternion polynomial of the form $g(q) = q^{n} +q^{n-1}a_{n-1}+...+qa_1 + a_0$ is given by

   $$C_g=\begin{bmatrix}
0& 0& 0& ... & 0& -a_0\\
1& 0& 0& ... & 0& -a_1\\
0& 1& 0& ... & 0& -a_2\\
.& .& .& ...& .& .\\
0& 0& 0& ...& 1& -a_{n-1}\\
\end{bmatrix}$$\\
Also Quaternion Companion Matrix for a monic lacunary type  polynomial of the form \\ $p(q) = q^{n} +a_{r}q^{r}+...+a_{1}q + a_{0}$,where $a_{r}\neq0$
and~~$1\leq{r}\leq{n-1}$  
is given by \\

	$$C_{p}=\begin{bmatrix}
0&1&0&...&0&.&0\\
0&0&1&...&0&.&0\\
0&0&0&...&0&.&0\\
.&.&.&.&.&.&.\\
.&.&.&.&.&.&.\\
.&.&.&.&.&.&.\\
0&0&0&...&0&.&1\\
-a_{0}& -a_{1}&-a_{2}&...& -a_{r}&.& 0\\
\end{bmatrix}$$

\parindent=0mm\vspace{0.in}Whereas the  companion matrix for a monic lacunary type  quaternion polynomial of the form  \\
$h(q) = q^{n} +q^{r}a_{r}+...+qa_1 + a_0$,
where $a_{r}\neq0$
and~~$1\leq{r}\leq{n-1}$  
is given by \\
 $$C_h=\begin{bmatrix}
0& 0& 0&...&0&...&0&-a_{0} \\
1& 0& 0&...&0&...&0&-a_{1}\\
.& .& .&...&.&...&.&.\\
0& 0& 0&...&1&...&0&-a_{r}\\
.& .& .&...&.&...&.&.\\
0& 0& 0&...&0&...&1&0\\
\end{bmatrix}$$

 For an $n \times n$ matrix $A=(p_{\mu\nu})$ of quaternions, the non commutativity of quaternions result in two
different types of eigenvalues,\cite{ee}.

\parindent=0mm\vspace{0.in}

{\bf{Quaternionic eigenvalues.}} Let $A$ be an $n \times n$ matrix with
entries from $\mathbb{H}$. Lack of commutativity means that we have two
possible types of eigenvalue viz left eigenvalue and right eigenvalue.

A quaternion $\lambda$ is said to be left eigenvalue of $A\in\mathbb{H}^{n\times n}[q]$, if $Ax=\lambda x$ for some non-zero  quaternionic vector $x \in \mathbb{H}^n$ and the vector $x$ is called a left eigenvector associated with it and a quaternion $\lambda$ is said to be right eigenvalue of $A\in\mathbb{H}^{n\times n}[q]$ if $Ax=x\lambda$ for some non-zero  quaternionic vector $x \in \mathbb{H}^n$ and the vector $x$ is called a right eigenvector associated with it. As usual the left spectrum and right spectrum of $A\in\mathbb{H}^{n\times n}[q]$ are respectively defined by
$$\sigma_{l}(A)=\{\lambda\in\mathbb{H}:Ax=\lambda x, for ~~some~~ x\neq0\}$$
and 
$$\sigma_{r}(A)=\{\lambda\in\mathbb{H}:Ax=x\lambda , for ~~some~~ x\neq0\}.$$
For complex case, concerning the location of eigenvalues, the famous Ger\v{s}gorin theorem can be stated as
\begin{thmx}\label{A}
All the eigenvalues of a $n \times n$ complex matrix $A=(a_{\mu\nu})$ are contained in the union of n Ger\v{s}gorin discs defined by $D_{\mu}=\{ z \in \mathbb{C}:|z-a_{\mu\mu}|\leq  \sum_{{\nu=1}\atop{\nu\neq\mu}}^n |a_{\mu\nu}|\}$.
\end{thmx}
Recently, Dar et al. \cite{RD}  extended Theorem \ref{A} to the quaternions, more precisely they proved following quaternion  version of Ger\v{s}gorin theorem.\\
\begin{thmx}\label{B}
All the left eigenvalues of a $n \times n$ matrix $A=(a_{\mu\nu})$ of quaternions lie  in the union of the n Ger\v{s}gorin balls defined by $B_{\mu}=\{q \in \mathbb{H}:|q-a_{\mu\mu}|\leq \rho_{\mu}(A) \},$ where $\rho_{\mu}(A)=\sum_{{\nu=1}\atop{\nu\neq\mu}}^n |a_{\mu\nu}|.$
\end{thmx}
Since the set of zeros of quaternion polynomial and the left spectrum of its companion matrix coincide, i.e, if $h(q)$ is a quaternion polynomial then $Zero(h)=\sigma_{l}(C_{h})$ \cite{gg}, we may therefore use Theorem \ref{B} and other known results on the left eigenvalues of quaternionic matrices as a tool in determining the zeros of a given polynomial and vice-versa.\\

We study analytic theory of polynomials in geometric theory of functions, where we focus on polynomials $P(z)=\sum\limits_{j=o}^{n}a_j z^j$ of a complex variable and the different parameters connected with it, such as coefficients, zeros, bounds for the zeros, maximum modulus of $P(z)$, location of the zeros, the relationship between the zeros and coefficients, zero free regions, growth, and much more. In accordance with the Fundamental Theorem of Algebra, Every polynomial of positive degree with complex coefficients always has a zero and the number of zeros is exactly equal to its degree. This theorem although tells about the exact number of zeros of a polynomial but the zeros can’t be determined
algebraically when the degree of a polynomial is more than 4. Thus, the challenge of determining at least the regions of a polynomial that include all or some of its zeros becomes important in the theory of polynomials.The first contributors of the subject were Gauss and Cauchy. As Compared to Gauss's bounds, Cauchy's were more exact for the moduli of the zeros of polynomials having arbitrary complex coefficients. In 1829,concerning the location of zeros of a polynomial with complex coefficients,  
A. L Cauchy \cite{ff} gave
a very simple expression of the region containing all the zeros in terms of the coefficients of a
polynomial. He proved the following theorem.
\begin{thmx}\label{C}
If $P(z)=z^n+z^{n-1}c_{n-1}+z^{n-2}c_{n-2}+...+zc_1+c_0$, is a complex polynomial of degree $n$, then all the zeroes of $P(z)$ lie inside the disc $|z|<1+\max\limits_{0\leq\nu\leq n-1}|c_\nu|$.
\end{thmx}
M.Fujiwara \cite{hh}, generalized Theorem \ref{C} by proving the following result.
\begin{thmx}\label{D}
If $ \lambda_1,\lambda_2,\lambda_3,...,\lambda_n $ are positive numbers such that $ \lambda_1+\lambda_2+\lambda_3+...+\lambda_n=1 $ .Then all the zeroes of the polynomial $P(z)=z^n+a_nz^{n-1}+a_{n-1}z^{n-2}+...+a_1z+a_0$ of degree $n$ lie in the circle.
$$|z|\leq Max\bigg(\frac{|a_{n-j}|}{\lambda_j}\bigg)^{\frac{1}{j}} ,~~~~~~~j=1,2,...,n$$
\end{thmx}
Recently, Dar et al. \cite{RD}  extended Theorem \ref{C} to the quaternions settings by proving the following result.
\begin{thmx}\label{E}
If $g(q)=q^n+q^{n-1}a_{n-1}+q^{n-2}a_{n-2}+...+qa_1+a_0$, is a quaternion polynomial with quaternion coefficients and $q$ is quaternionic variable, then all the zeroes of $g(q)$ lie in the ball $|q|<1+\max\limits_{0\leq\nu\leq n-1}|a_\nu|$.
\end{thmx}
\section{\textbf{Main Results}}
In this paper, we will use various matrix tools including Theorem \ref{B} to prove the various results concerning the location of zeros of quaternionic polynomials. We begin by proving the following generalization of Theorem \ref{E}.
\begin{theorem}\label{TH1}
 Let $g(q) = q^{n} +q^{n-1}a_{n-1}+...+qa_1 + a_0$ be a quaternion polynomial with quaternion coefficients and $q$ be a quaternion variable, then for any positive numbers $\alpha_1,\alpha_2,...,\alpha_{n-1}$ with $\alpha_0=0$ and $\alpha_n=1,$ all zeroes of $g(q)$ lie in the ball

\begin{align*}
\bigg\{ q \in \mathbb{H} : |q| \leq \max \bigg\{ \frac{\alpha_i}{\alpha_{i+1}} + \frac{|a_i|}{\alpha_{i+1}}: \ i=0,1,2,...,n-1 \bigg\} \bigg\}
\end{align*}

\end{theorem}
\begin{remark}\label{R1}
For $\alpha_1=\alpha_2=\alpha_3= ... =\alpha_n = 1,$ Theorem \ref{TH1} reduces to Theorem \ref{E}. On the other hand if we take $$\alpha_1=|a_1|, \alpha_2=|a_2|,...,\alpha_{n-1}=|a_{n-1}|.$$
Then
$$\frac{\alpha_1}{\alpha_2}=\frac{|a_1|}{|a_{2}|},\frac{\alpha_2}{\alpha_3}=\frac{|a_2|}{|a_{3}|},...,\frac{\alpha_{n-1}}{\alpha_n}=|a_{n-1}|.$$\\
In the view of this and  the fact that $\alpha_{n}=1$, Theorem \ref{TH1} reduces to the following result.
\end{remark}
\begin{corollary}\label{C1}
 All the zeroes of polynomial $g(q) = q^{n} +q^{n-1}a_{n-1}+...+qa_1 + a_0$ lie in the ball 
$$|q| \leq \max\bigg\{ \bigg|\frac{a_0}{a_1} \bigg|,2 \bigg|\frac{a_1}{a_2} \bigg|,...,2|a_{n-1}|\bigg\}.$$
\end{corollary}
Next, we extend Theorem \ref{D}  to quaternionic settings, More precisely we prove:
\begin{theorem}\label{TH2}
If $ \lambda_1,\lambda_2,\lambda_3,...,\lambda_n $ are positive numbers such that $ \lambda_1+\lambda_2+\lambda_3+...+\lambda_n=1 $ .Then all the zeros of the quaternionic polynomial $f(q)=q^n+a_{n-1}q^{n-1}+a_{n-2}q^{n-2}+...+a_1q+a_0$ lie in the ball 
$$\bigg\{q\in \mathbb{H} :|q|\leq \max \bigg\{\bigg(\frac{|a_{n-j}|}{\lambda_j}\bigg)^{\frac{1}{j}} :~~~~j=1,2,...,n\bigg\}\bigg\}.$$
\end{theorem}
\indent Now we turn towards lucunary type of quaternionic polynomials and prove some interesting results, which in turn will yield various extensions of Theorem \ref{C} to the quaternionic settings including Theorem \ref{TH1} as special cases.
\begin{theorem}\label{TH3}
Let $p(q)=q^n+a_rq^r+...+a_1q+a_0$, $a_r\neq 0,~ 0\leq r\leq n-1$ be a quaternionic polynomial and $q$ be the quaternion variable. If $\lambda=\bigg\{\max\limits_{0\leq j\ \leq r}|a_j|\bigg\}^{\frac{1}{n}},$ then all the zeros of $p(q)$ lie in the ball
\begin{align}\label{Te1}
\big\{q\in \mathbb{H}:|q|\leq (\lambda+\lambda^2+...+\lambda^{r+1}),where ~~1\leq r\leq n-1\big \}
\end{align}
\end{theorem}
\begin{theorem}\label{TH4}
Let $h(q)=q^n+q^ra_r+...+qa_1+a_0$, $a_r\neq 0,~ 0\leq r\leq n-1$ be a quaternionic polynomial and $q$ be the quaternion variable. If $\lambda=\bigg\{\max\limits_{0\leq j\ \leq r}|a_j|\bigg\}^{\frac{1}{n}},$ then all the zeros of $h(q)$ lie in the ball
 \begin{align}\label{Te2}
\big\{q\in \mathbb{H}:|q|\leq \lambda+\max(\lambda^2,\lambda^{r+1}),~~where ~~1\leq r\leq n-1 \big \} 
 \end{align}
\end{theorem}
\begin{remark}\label{R2}
It is easy to see that for $1 \leq r< n-1,$ the bound given by \eqref{Te2} is stronger than \eqref{Te1}. On the other hand if we take $r=n-2$ in Theorem \ref{TH3} and Theorem \ref{TH4}, we get following two results.
\end{remark}
\begin{corollary}
All the zeros of the quaternionic polynomial $p(q)=q^n+a_{n-2}q^{n-2}+...+a_1q+a_0$ lie in the ball
$$|q|\leq\lambda+\lambda^{2}+\lambda^{3}+...+\lambda^{n-1},$$ where $\lambda $ is defined in Theorem \ref{TH3}.
\end{corollary}
\begin{corollary}
All the zeros of the quaternionic polynomial $h(q)=q^n+q^{n-2}a_{n-2}+...+qa_1+a_0$ lie in the ball
$$|q|\leq\lambda+\max\{\lambda^{2},\lambda^{n-1}\},$$ where $\lambda $ is defined in Theorem \ref{TH3}.
\end{corollary}
\section{\textbf{Auxiliary Results}}
For the proof of our main results, we need the following two Lemmas of Dar et al. \cite{RD} and Rather et al\cite{R}.
\begin{lemma}\label{L1}
All the left eigenvalues of a $n \times n$ matrix $A=(a_{\mu\nu})$ of quaternions lie in the union of the n Ger\v{s}gorin balls defined by $B_{\mu}=\{q \in \mathbb{H}:|q-a_{\mu\mu}|\leq \rho_{\mu}(A) \}$ where $\rho_{\mu}(A)=\sum_{{\nu=1}\atop{\nu\neq\mu}}^n |a_{\mu\nu}|.$
\end{lemma}
\begin{lemma}\label{L2}
Let $P(q)$ be a quaternion polynomial with quaternionic coefficients and $C$ be the companion matrix of $P(q)$, then for any diagonal matrix $D= diag(d_1, d_2, ..., d_{n-1}, d_n),$ where $d_1, d_2, ..., d_n$ are positive real numbers, the left eigenvalues of $D^{-1}C D$ and the zeros of $P(q)$ are same.
\end{lemma}

\section{\textbf{Proof of the main theorems}}
\textbf{Proof of theorem \ref{TH1}: }  Let $C_g$ be the companion matrix of the polynomial $g(q)$ and $T=diag[\alpha_1,\alpha_2,...,\alpha_{n-1}, 1]$ be a diagonal matrix.
Then 
$$T^{-1}C_gT=\begin{bmatrix} 
\frac{1}{\alpha_1} & 0 & 0 & 0 & . & . & . & 0 & 0  \\
0 & \frac{1}{\alpha_2} & 0 & 0 & . & . & . & 0 & 0\\. & . & . & . & . & . & . & . & .\\ 0 & 0 & 0 & 0 & . & . & . &\frac{1}{\alpha_{n-1}} & 0 \\
0 & 0 & 0 & 0 & . & . & . & . &  1 \\   
	\end{bmatrix}
	\begin{bmatrix}
0& 0& 0& ... & 0& -a_0\\
1& 0& 0& ... & 0& -a_{1}\\
0& 1& 0& ... & 0& -a_{2}\\
.& .& .& ...& .& .\\
0& 0& 0& ...& 1& -a_{n-1}\\
\end{bmatrix}
\begin{bmatrix} 
\alpha_1 & 0 & 0 & 0 & . & . & . & 0 & 0  \\
0 & \alpha_2 & 0 & 0 & . & . & . & 0 & 0\\. & . & . & . & . & . & . & . & .\\ 0 & 0 & 0 & 0 & . & . & . &\alpha_{n-1} & 0 \\
0 & 0 & 0 & 0 & . & . & . & . &  1 \\   
	\end{bmatrix}$$

$$=\begin{bmatrix}
0& 0& 0& ...& 0& \frac{-a_{0}}{\alpha_1}\\
\frac{\alpha_1}{\alpha_2}& 0& 0& ...& 0& \frac{-a_1}{\alpha_2}\\
0& \frac{\alpha_2}{\alpha_3}& 0& ...& 0& \frac{-a_2}{\alpha_3}\\
.& .& .& ...& .& .\\
0& 0& 0& ...&       \alpha_{n-1}&     -a_{n-1}\\
\end{bmatrix}$$
Applying Lemma \ref{L1} to the matrix $T^{-1}C_gT$, it follows that all the left eigenvalues of the matrix $T^{-1}C_gT$ lie in the union of balls.
$$|q|\leq\frac{|a_0|}{\alpha_1}<\frac{\alpha_0}{\alpha_1}+\frac{|a_0|}{\alpha_1}$$
$$|q|\leq\frac{\alpha_1}{\alpha_2}+\frac{|a_1|}{\alpha_2}$$
$$------$$
$$------$$and
$$|q|\leq \alpha_{n-1}+|a_{n-1}|.$$
That is, all the left eigenvalues of the matrix  $T^{-1}C_gT$ lie in the union of the balls
\begin{align}\label{te11}
|q| \leq \left\{ \frac{\alpha_i}{\alpha_{i+1}} + \frac{|a_i|}{\alpha_{i+1}} \right\}, \quad i = 0,1,2, \ldots, n-1~ \textnormal{with}~ \alpha_n=1.
\end{align}
In other words, all the left eigenvalues of the matrix $T^{-1}C_gT$ lie in the ball
\begin{align}\label{k1}
|q| \leq \max \left\{ \frac{\alpha_i}{\alpha_{i+1}} + \frac{|a_i|}{\alpha_{i+1}}, i = 0,1,2, \ldots, n-1\right\} ~ \textnormal{where}~ \alpha_n=1.
\end{align}

Since $T$ is a diagonal matrix with real positive entries, by Lemma \ref{L2}, it follows that the left
eigenvalues of $T^{-1}C_gT$ are the zeros of $g(q).$ Hence, all the zeros of $g(q)$ lie in the ball given by \eqref{k1}.\\
That completes the proof of Theorem \ref{TH1}.
 
\textbf{Proof of theorem \ref{TH2}: } Let $C_f$ be the companion matrix of the polynomial $f(q)$ and $T=diag \left[\frac{1}{l^{n-1}},\frac{1}{l^{n-2}},...,\frac{1}{l},1\right]$ be a diagonal matrix , where $l$ is a positive real number.
Then 
$$T^{-1}C_fT=\begin{bmatrix} 
0 & l & 0 & 0 & . & . & . & 0 & 0  \\
0 & 0 & l & 0 & . & . & . & 0 & 0\\0 & 0 & 0 & l& . & . & . & 0 & 0\\. & . & . & . & . & . & . &. & .\\ 0 & 0 & 0 & 0 & . & . & . &0 &l\\
-\frac{a_{0}}{l^{n-1}} & -\frac{a_1}{l^{n-2}} & -\frac{a_2}{l^{n-3}} & -\frac{a_3}{l^{n-4}} & . & . & . & -\frac{a_{n-2}}{l} & -a_{n-1}\\   
	\end{bmatrix}$$

Applying Lemma \ref{L1} to the matrix $T^{-1}C_fT$, it follows that all the left eigenvalues of the matrix $T^{-1}C_fT$ lie in the union of balls
$|q|\leq l$ and $|q+a_{n-1}|\leq \frac{|a_0|}{l^{n-1}}+ \frac{|a_1|}{l^{n-2}}+...+\frac{|a_{n-2}|}{l}.$ 

Since 
\begin{align*}
\begin{split}
|q|{}&= |q+a_{n-1}-a_{n-1}|\\
&\leq|s+a_{n-1}|+|a_{n-1}|\\
&\leq \sum_{j=1}^n\frac{|a_{n-j}|}{l^{j-1}}.
\end{split}
\end{align*}
That is, all the left eigenvalues of the matrix $T^{-1}C_fT$  lie in the ball
\begin{align}\label{t2e1}
|q|\leq \max \big\{l,\sum_{j=1}^n\frac{|a_{n-j}|}{l^{j-1}}\big\}.
\end{align}
We now choose
$$l= \max \bigg\{\frac{|a_{n-j}|}{\lambda_{j}}\bigg\}^{\frac{1}{j}}, ~~~~j=1,2,3,..,n$$
then $$\frac{|{a_{n-j}}|}{\lambda_j}\leq l^j \quad \forall \quad j=1,2,3,...n,$$
which gives
$$\frac{|a_{n-j}|}{l^{j-1}}\leq l\cdot\lambda_j$$

so that

$$\sum_{j=1}^n\frac{|a_{n-j}|}{\l^{j-1}} \leq  \sum_{j-1}^nl\lambda_j=l(\lambda_1+\lambda_2+\lambda_3+...+\lambda_n)=l.$$   
Using this in \eqref{t2e1}, it follows that all the left eigenvalues of $T^{-1}C_{f}T$ lie in the ball
\begin{align}\label{t2e2}
|q|\leq \max \bigg\{ \bigg( \frac{|a_{n-j}|}{\lambda_j}\bigg)^{\frac{1}{j}},~ j=1,2,...,n\bigg\}.
\end{align}
Since $T$ is a diagonal matrix with real positive entries, by Lemma \ref{L2}, it follows that the left
eigenvalues of $T^{-1}C_fT$ are the zeros of $f(q).$ Hence, all the zeros of $f(q)$ lie in the ball given by \eqref{t2e2}.\\
That completes the proof of Theorem \ref{TH2}.\\
 {\bf{Proof of Theorem \ref{TH3}:}} ~~Let $C_{p}$ be the companion matrix of the polynomial $p(q),$ we take a matrix  $$T=diag \left[ \frac{1}{\lambda^{n-1}},\frac{1}{\lambda^{n-2}},...,\frac{1}{\lambda},1\right]$$
 where $\lambda$ is positive real number 
  then 
$$T^{-1}C_pT=\begin{bmatrix} 
\lambda^{n-1} & 0 & 0 & 0 & . & . & . & 0 & 0  \\
0 & \lambda^{n-2} & 0 & 0 & . & . & . & 0 & 0\\. & . & . & . & . & . & . & . & .\\
. & . & . & . & . & . & . & . & .\\
. & . & . & . & . & . & . & . & .\\ 0 & 0 & 0 & 0 & . & . & . &\lambda & 0 \\
0 & 0 & 0 & 0 & . & . & . & . &  1 \\   
	\end{bmatrix}
	\begin{bmatrix}
0&1&0&...&0&.&0\\
0&0&1&...&0&.&0\\
0&0&0&...&0&.&0\\
.&.&.&.&.&.&.\\
.&.&.&.&.&.&.\\
.&.&.&.&.&.&.\\
0&0&0&...&0&.&1\\
-a_{0}& -a_{1}&-a_{2}&...& -a_{r}&.& 0\\
\end{bmatrix}
\begin{bmatrix} 
\frac{1}{\lambda^{n-1}} & 0 & 0 & 0 & . & . & . & 0 & 0  \\
0 & \frac{1}{\lambda^{n-2}} & 0 & 0 & . & . & . & 0 & 0\\. & . & . & . & . & . & . & . & .\\
. & . & . & . & . & . & . & . & .\\
. & . & . & . & . & . & . & . & .\\ 0 & 0 & 0 & 0 & . & . & . &\frac{1}{\lambda} & 0 \\
0 & 0 & 0 & 0 & . & . & . & . &  1 \\   
\end{bmatrix}
$$
	
	$$=\begin{bmatrix}
0& \lambda & ...& 0&...& 0\\
0& 0 & ...& 0&...& 0\\
.& . & ...& .&...& .\\
0& 0 & ...& \lambda&...& 0\\
\frac{-a_0}{\lambda^{n-1}}& \frac{-a_1}{\lambda^{n-2}} & ...& \frac{-a_r}{\lambda^{n-r-1}}&...& 0\\
\end{bmatrix}$$
Applying Lemma \ref{L2} to the matrix $T^{-1}C_pT$. It follows all the left eigenvalues of the matrix $T^{-1}C_pT$ lie in the union of the balls $|q|\leq \lambda$ and $$|q|\leq \bigg |\frac{a_0}{\lambda^{n-1}}\bigg |+\bigg |\frac{a_1}{\lambda^{n-2}}\bigg |+...+\bigg |\frac{a_r}{\lambda^{n-r-1}}\bigg |=\sum_{j=0}^r\frac{|a_j|}{\lambda^{n-j-1}}$$
that is , all the left eigenvalues of the matrix $T^{-1}C_pT$ lie in the ball
\begin{align}\label{t3e1}
|q|\leq \max \bigg \{\lambda,\sum_{j=0}^r\frac{|a_j|}{\lambda^{n-j-1}}\bigg \}.
\end{align}
Since by hypothesis
$$\lambda^n=\max\limits_{0\leq j\ \leq r}|a_j|, ~~~\lambda \neq 0$$
we have $$|a_j| \leq \lambda^n.$$
Therefore from \eqref{t3e1}, we have
\begin{eqnarray*}
|q|&\leq & \max \bigg \{\lambda,\sum_{j=0}^r\lambda^{j+1}\bigg \}\\
&=& \max \{\lambda,\lambda+\lambda^2+\lambda^3+...\lambda^{r+1}\}\\
&=& \lambda+\lambda^2+...+\lambda^{r+1},
\end{eqnarray*} 
that is,
\begin{align}\label{te22} 
|q| \leq \lambda+\lambda^2+...+\lambda^{r+1}.
\end{align}
Since $T$ is a diagonal matrix with real positive entries, by Lemma \ref{L2}, it follows that the left
eigenvalues of $T^{-1}C_pT$ are the zeros of $p(q).$ Hence, all the zeros of $p(q)$ lie in the ball given by \eqref{te22}.\\
That completes the proof of Theorem \ref{TH3}.\\
{\bf{Proof of Theorem \ref{TH4}:}}\\
Let $C_h$ be the companion matrix of $h(q) = q^{n}a_n +q^{r}a_{r}+...+qa_1 + a_0, a_r\neq 0, 0\leq r\leq n-1$ and let $T=diag \left[ \lambda^{n-1},\lambda^{n-2},...,\lambda,1\right]$ be a diagonal matrix with positive entries, then  
$$T^{-1}C_{h}T =\begin{bmatrix}
0& 0& ... & 0 &...&0& \frac{-a_{0}}{\lambda^{n-1}}\\
\lambda& 0& ... & 0 &... &0&\frac{-a_{1}}{\lambda^{n-2}}\\
.& .& ... & . &... &.&.\\
0& 0& ... & \lambda &... &0&\frac{-a_{r}}{\lambda^{n-r-1}}\\
.& .& ... & . &...&.& .\\
0&0& ... & 0 &...&\lambda& 0\\
\end{bmatrix}$$
Applying Lemma \ref{L1} to the matrix $T^{-1}C_hT,$ it follows that all the left eigenvalues of the matrix $T^{-1}C_hT$ lie in the ball
$$|q|\leq\max\limits_{1\leq\j\leq r}\bigg\{\frac{|a_{0}|}{\lambda^{n-1}},\lambda+\frac{|a_{j}|}{\lambda^{n-j-1}}\bigg\}$$\\
     $$\leq\max\limits_{1\leq\j\leq r}\{\lambda,\lambda+\lambda^{j+1}\},$$
that is, all the left eigenvalues of the matrix $T^{-1}C_hT$ lie in the ball
\begin{align*}
\begin{split}
|q|{}& \leq \max\{ \lambda, \lambda + \lambda^{2}, \lambda + \lambda^{3}, ..., \lambda + \lambda^{r+1} \} \\
&= \lambda + \max\{ \lambda^{2}, \lambda^{3}, ..., \lambda^{r+1} \}\\
& = \lambda + \max(\lambda^2, \lambda^{r+1}).
\end{split}
\end{align*}

Since $T$ is a diagonal matrix with real positive entries, by Lemma \ref{L2}, it follows that the left
eigenvalues of $T^{-1}C_hT$ are the zeros of $h(q).$ Hence, all the zeros of $h(q)$ lie in the ball given by \eqref{Te2}.\\
That completes the proof of Theorem \ref{TH4}.\\

\end{document}